\documentclass{amsart}
\usepackage{latexsym}
\usepackage{amssymb}
\usepackage{amsmath}

\begin{document}


\newtheorem{theorem}{Theorem}
\newtheorem{problem}{Problem}
\newtheorem{definition}{Definition}
\newtheorem{lemma}{Lemma}
\newtheorem{proposition}{Proposition}
\newtheorem{corollary}{Corollary}
\newtheorem{example}{Example}
\newtheorem{conjecture}{Conjecture}
\newtheorem{algorithm}{Algorithm}
\newtheorem{exercise}{Exercise}
\newtheorem{remarkk}{Remark}

\newcommand{\be}{\begin{equation}}
\newcommand{\ee}{\end{equation}}
\newcommand{\bea}{\begin{eqnarray}}
\newcommand{\eea}{\end{eqnarray}}
\newcommand{\beq}[1]{\begin{equation}\label{#1}}
\newcommand{\eeq}{\end{equation}}
\newcommand{\beqn}[1]{\begin{eqnarray}\label{#1}}
\newcommand{\eeqn}{\end{eqnarray}}
\newcommand{\beaa}{\begin{eqnarray*}}
\newcommand{\eeaa}{\end{eqnarray*}}
\newcommand{\req}[1]{(\ref{#1})}

\newcommand{\lip}{\langle}
\newcommand{\rip}{\rangle}

\newcommand{\uu}{\underline}
\newcommand{\oo}{\overline}
\newcommand{\La}{\Lambda}
\newcommand{\la}{\lambda}
\newcommand{\eps}{\varepsilon}
\newcommand{\om}{\omega}
\newcommand{\Om}{\Omega}
\newcommand{\ga}{\gamma}
\newcommand{\rrr}{{\Bigr)}}
\newcommand{\qqq}{{\Bigl\|}}

\newcommand{\dint}{\displaystyle\int}
\newcommand{\dsum}{\displaystyle\sum}
\newcommand{\dfr}{\displaystyle\frac}
\newcommand{\bige}{\mbox{\Large\it e}}
\newcommand{\integers}{{\Bbb Z}}
\newcommand{\rationals}{{\Bbb Q}}
\newcommand{\reals}{{\rm I\!R}}
\newcommand{\realsd}{\reals^d}
\newcommand{\realsn}{\reals^n}
\newcommand{\NN}{{\rm I\!N}}
\newcommand{\degree}{{\scriptscriptstyle \circ }}
\newcommand{\dfn}{\stackrel{\triangle}{=}}
\def\complex{\mathop{\raise .45ex\hbox{${\bf\scriptstyle{|}}$}
     \kern -0.40em {\rm \textstyle{C}}}\nolimits}
\def\hilbert{\mathop{\raise .21ex\hbox{$\bigcirc$}}\kern -1.005em {\rm\textstyle{H}}} 
\newcommand{\RAISE}{{\:\raisebox{.6ex}{$\scriptstyle{>}$}\raisebox{-.3ex}
           {$\scriptstyle{\!\!\!\!\!<}\:$}}} 

\newcommand{\hh}{{\:\raisebox{1.8ex}{$\scriptstyle{\degree}$}\raisebox{.0ex}
           {$\textstyle{\!\!\!\! H}$}}}

\newcommand{\OO}{\won}
\newcommand{\calA}{{\mathcal A}}
\newcommand{\BB}{{\mathcal B}}
\newcommand{\calC}{{\cal C}}
\newcommand{\calD}{{\cal D}}
\newcommand{\calE}{{\cal E}}
\newcommand{\calF}{{\mathcal F}}
\newcommand{\calG}{{\cal G}}
\newcommand{\calH}{{\cal H}}
\newcommand{\calK}{{\cal K}}
\newcommand{\calL}{{\mathcal L}}
\newcommand{\calM}{{\cal M}}
\newcommand{\calO}{{\cal O}}
\newcommand{\calP}{{\cal P}}
\newcommand{\calU}{{\mathcal U}}
\newcommand{\calX}{{\cal X}}
\newcommand{\calXX}{{\cal X\mbox{\raisebox{.3ex}{$\!\!\!\!\!-$}}}}
\newcommand{\calXXX}{{\cal X\!\!\!\!\!-}}
\newcommand{\gi}{{\raisebox{.0ex}{$\scriptscriptstyle{\cal X}$}
\raisebox{.1ex} {$\scriptstyle{\!\!\!\!-}\:$}}}
\newcommand{\intsim}{\int_0^1\!\!\!\!\!\!\!\!\!\sim}
\newcommand{\intsimt}{\int_0^t\!\!\!\!\!\!\!\!\!\sim}
\newcommand{\pp}{{\partial}}
\newcommand{\al}{{\alpha}}
\newcommand{\sB}{{\cal B}}
\newcommand{\sL}{{\cal L}}
\newcommand{\sF}{{\cal F}}
\newcommand{\sE}{{\cal E}}
\newcommand{\sX}{{\cal X}}
\newcommand{\R}{{\rm I\!R}}
\renewcommand{\L}{{\rm I\!L}}
\newcommand{\vp}{\varphi}
\newcommand{\N}{{\rm I\!N}}
\def\ooo{\lip}
\def\ccc{\rip}
\newcommand{\ot}{\hat\otimes}
\newcommand{\rP}{{\Bbb P}}
\newcommand{\bfcdot}{{\mbox{\boldmath$\cdot$}}}

\renewcommand{\varrho}{{\ell}}
\newcommand{\dett}{{\textstyle{\det_2}}}
\newcommand{\sign}{{\mbox{\rm sign}}}
\newcommand{\TE}{{\rm TE}}
\newcommand{\TA}{{\rm TA}}
\newcommand{\E}{{\rm E\,}}
\newcommand{\won}{{\mbox{\bf 1}}}
\newcommand{\Lebn}{{\rm Leb}_n}
\newcommand{\Prob}{{\rm Prob\,}}
\newcommand{\sinc}{{\rm sinc\,}}
\newcommand{\ctg}{{\rm ctg\,}}
\newcommand{\loc}{{\rm loc}}
\newcommand{\trace}{{\,\,\rm trace\,\,}}
\newcommand{\Dom}{{\rm Dom}}
\newcommand{\ifff}{\mbox{\ if and only if\ }}
\newcommand{\nproof}{\noindent {\bf Proof:\ }}
\newcommand{\remark}{\noindent {\bf Remark:\ }}
\newcommand{\remarks}{\noindent {\bf Remarks:\ }}
\newcommand{\note}{\noindent {\bf Note:\ }}

\newcommand{\boldx}{{\bf x}}
\newcommand{\boldX}{{\bf X}}
\newcommand{\boldy}{{\bf y}}
\newcommand{\boldR}{{\bf R}}
\newcommand{\uux}{\uu{x}}
\newcommand{\uuY}{\uu{Y}}

\newcommand{\limn}{\lim_{n \rightarrow \infty}}
\newcommand{\limN}{\lim_{N \rightarrow \infty}}
\newcommand{\limr}{\lim_{r \rightarrow \infty}}
\newcommand{\limd}{\lim_{\delta \rightarrow \infty}}
\newcommand{\limM}{\lim_{M \rightarrow \infty}}
\newcommand{\limsupn}{\limsup_{n \rightarrow \infty}}

\newcommand{\ra}{ \rightarrow }

\newcommand{\ARROW}[1]
  {\begin{array}[t]{c}  \longrightarrow \\[-0.2cm] \textstyle{#1} \end{array} }

\newcommand{\AR}
 {\begin{array}[t]{c}
  \longrightarrow \\[-0.3cm]
  \scriptstyle {n\rightarrow \infty}
  \end{array}}

\newcommand{\pile}[2]
  {\left( \begin{array}{c}  {#1}\\[-0.2cm] {#2} \end{array} \right) }

\newcommand{\floor}[1]{\left\lfloor #1 \right\rfloor}

\newcommand{\mmbox}[1]{\mbox{\scriptsize{#1}}}

\newcommand{\ffrac}[2]
  {\left( \frac{#1}{#2} \right)}

\newcommand{\one}{\frac{1}{n}\:}
\newcommand{\half}{\frac{1}{2}\:}

\def\le{\leq}
\def\ge{\geq}
\def\lt{<}
\def\gt{>}

\def\squarebox#1{\hbox to #1{\hfill\vbox to #1{\vfill}}}
\newcommand{\nqed}{\hspace*{\fill}
           \vbox{\hrule\hbox{\vrule\squarebox{.667em}\vrule}\hrule}\bigskip}

\newcommand{\no}{\noindent}
\newcommand{\EE}{\mathbb{E}}
\newcommand{\RR}{\mathbb{R}}
\newcommand{\D}{\mathcal{D}}
\newcommand{\DD}{\widetilde{\mathcal{D}}}
\newcommand{\DDD}{\widetilde{\widetilde{\mathcal{D}}}}
\newcommand{\LL}{\mathcal{L}}
\newcommand{\F}{\mathcal{F}}
\newcommand{\B}{\mathcal{B}}
\newcommand{\M}{\mathcal{M}}
\newcommand{\W}{\mathbb{W}}
\newcommand{\lo}{lp}
\newcommand{\grandX}{\mathbb{X}}
\newcommand{\loigrandX}{\mu^{\mathbb{X}}}
\newcommand{\mup}{\mu_*}
\newcommand{\betap}{\beta_*}
\newcommand{\pa}{pa}
\newcommand{\QQ}{\mathbb{Q}}


\title{VARIATIONAL CALCULUS ON WIENER SPACE WITH RESPECT TO CONDITIONAL EXPECTATIONS}

\author{ K\'evin HARTMANN}
\maketitle
\noindent
{\bf Abstract:}{\small{We give a variational formulation for $-\log\EE_\nu\left[e^{-f}|\F_t\right]$ for a large class of measures $\nu$. We give a refined entropic characterization of the invertibility of some perturbations of the identity. We also discuss the attainability of the infimum in the variational formulation and obtain a Pr\'ekopa-Leindler theorem for conditional expectations.
}}\\

\vspace{0.5cm}


\tableofcontents

\section{\bf{Introduction}}

Denote $\W$ the space of continuous functions from $[0,1]$ to $\RR^n$ and H the associated canonical Cameron-Martin space of elements of W which admit a density in $L^2$. Also denote $\mu$ the Wiener measure, W the coordinate process, and $(\F_t)$ the canonical filtration of W completed with respect to $\mu$. W is a Brownian motion under $\mu$. Set f a bounded from above measurable function from $\W$ to $\RR$. In \cite{du}, Dupuis and Ellis prove that
\bea -\log\EE_\mu\left[e^{-f}\right]=\inf_\theta\left(\EE_\theta\left[f\right]+H(\theta|\mu)\right)\label{3rv0}\eea
\no where the infimum is taken over the probability measures $\theta$ on $\W$ which are absolutely continuous with respect to $\mu$ and the relative entropy $H(\theta|\mu)$ is equal to $\EE_\mu\left[\frac{d\theta}{d\mu}\log\frac{d\theta}{d\mu}\right]$. In \cite{bd}, Bou\'e and Dupuis use it to derive the variational formulation
\bea -\log\EE_\mu\left[e^{-f}\right]=\inf_u\EE_\mu\left[f\circ (W+u)+\frac{1}{2}\int_0^1|\dot{u}(s)|^2ds\right]\label{3rv}\eea
\no where the infimum is taken over $L^2$ functions from $\W$ to H whose density is adapted to $(\F_t)$. This variational formulation is useful to derive large deviation asymptotics as Laplace principles for small noise diffusions for instance. This result was later extended by Budhiraja and Dupuis to Hilbert-space-valued Brownian motions in \cite{bud}, and then by Zhang to abstract Wiener space in \cite{zh}, using the framework developed by \"Ust\"unel and Zakai in \cite{ust1}.\newline
\no The Pr\'ekopa-Leindler theorem first formulation was given by Pr\'ekopa in \cite{pre1} and arose in stochastic programming where a lot of non-linear optimization problems require concavity. In \cite{ha}, Huu Hariya uses the variational formulation to retrieve a Pr\'ekopa-Leindler theorem for Wiener space, similar to the formulation of \"Ust\"unel in \cite{ust3} with log-concave measures. Other functional inequalities can be derived from \ref{3rv}, see for instance Lehec in \cite{le1}.\newline
\no The bounded from above hypothesis in \ref{3rv} was weakened significantly by \"Ust\"unel in \cite{art}, it was replaced with the condition
$$\EE_\mu\left[fe^{-f}\right]<\infty$$
\no and the existence of conjugate integers p and q such that
$$f\in L^p(\mu),e^{-f}\in L^q(\mu)$$
\no These relaxed hypothesis pave the way to new applications. The possibility of using unbounded functions is primordial in Dabrowski's application of \ref{3rv} to free entropy in \cite{da}.\newline
\no \"Ust\"unel's approach is routed in the study of the perturbations of the identity of $\W$, which is the coordinate process, and their invertibility. The question of the invertibility of an adapted perturbation of the identity is linked to the representability of measures and was put to light by the celebrated example of Tsirelson \cite{tsi}. \"Ust\"unel proved that if $u\in L^2(\mu,H)$ and its density is adapted, $I_\W+u$ is $\mu$-a.s. invertible if and only if
$$H((I_\W+u)\mu|\mu)=\frac{1}{2}\EE_\mu\left[|u|_H^2\right]$$
\no To prove \ref{3rv} with the integrability conditions specified above, \"Ust\"unel uses the fact that H-$C^1$ shifts, meaning shifts that are a.s. Fr\'echet-differentiable on H with a $\mu$-a.s. continuous on H Fr\'echet derivative, are a.s. invertible, and that shifts can be approached with H-$C^1$ shifts using the Ornstein-Uhlenbeck semigroup.\newline
\no In \cite{a1} we give a variational formulation similar as \ref{3rv} for diffusions solutions of stochastic differential equations, while lowering the integrability hypothesis on f.\newline
\no In \cite{a2} we present a general framework to be able to similarly derive a variational formulation for $-\log\EE_\nu[e^{-f}]$ for a large class of measures $\nu$, without increasing the integrability hypothesis on f. We give a set of conditions so that a set of processes $(W^u)$ can act as perturbations of W and allow a Girsanov-like change of variable with respect to a Brownian motion $\beta$. We write $\frac{e^{-f}}{\EE\left[e^{-f}\right]}$ as the Wick exponential of some v, and then approach v to obtain invertible perturbations of the identity.
\no Hyndman and Wang proved in \cite{hy} that
\bea-\log\EE_\mu\left[\left.e^{-f}\right|\F_t\right]=\inf_\theta\left(\EE_\theta\left[f|\F_t\right]+\EE_\theta\left[\left.\log\frac{d\theta}{d\mu}\right|\F_t\right]\right)\label{3rv2}\eea
\no where the infimum is taken over the probability measures $\theta$ which are absolutely continuous with respect to $\mu$ and verify $\EE_\mu\left[\frac{d\theta}{d\mu}|\F_t\right]=1$. They link it to forward-backward stochastic differential equations and apply it to various pricing problems for zero-coupon bonds.\newline
The relation \ref{3rv2}, obtained for a deterministic time t, is very similar to \ref{3rv0} so three questions arise naturally: can we obtain a relation similar to \ref{3rv} for the conditional expectation, can we extend it to other measures with the framework we developed in our third paper, and finally are these relations still valid if we substitute t with a stopping time $\tau$? Our paper answers affirmatively to these three questions. We keep the notations from our \cite{a2} and we prove that
\bea -\log\EE_\nu\left[\left.e^{-f}\right|\F_\tau\right]&=&\inf_\theta\left(\EE_\nu\left[f|\F_\tau\right]+\EE_\nu\left[\left.\frac{d\theta}{d\nu}\log\frac{d\theta}{d\nu}\right|\F_\tau\right]\right)\label{3rv3}\\
-\log\EE_\nu\left[\left.e^{-f}\right|\F_\tau\right]&=&\inf_u\EE_\nu\left[\left.f\circ W^u+\frac{1}{2}|u|_H^2\right|\F_\tau\right]\label{3rva}\eea
\no In \ref{3rv3} we assume that $\EE_\nu\left[fe^{-f}\right]<\infty$ and the infimum is taken over the probability measures $\theta$ on $\W$ which are absolutely continuous with respect to $\nu$ and such that $\EE_\nu\left[\frac{d\theta}{d\nu}|\F_t\right]=1$. In \ref{3rva}, the infimum is taken over the u from W to H, with adapted density, which are in $L^2$ and such that $1_{t\leq\tau}\dot{u}(t)=0$, and we assume that $\EE_\nu\left[fe^{-f}\right]<\infty$ and that there exists two conjugate integers p and q such that $f\in L^p(\nu)$ and $e^{-f}\in L^q(\nu)$. Observe that we had to increase the integrability hypothesis on f from what we had for the non-conditional case. In fact the integrability hypothesis on f here are the same as in \cite{art}. Finally, we discuss the attainability of the infimum in \ref{3rva} and we obtain a analog of Pr\'ekopa-Leindler type theorem for the conditional expectation with respect to $\mu$. However, similarly as in \cite{a2}, the convexity hypothesis seem quite restrictive.

\section{\bf{Framework}\label{fr3}}

\no Set $n\in\NN^*$, we denote $\W=C([0,1],\RR^{n})$ the canonical
Wiener space, $H=\left\{\int_0^.\dot{h}(s)ds, \dot{h}\in
  L^2([0,1])\right\}$ the associated Cameron-Martin space and W
is the coordinate process. We denote $(\F_t)$ its filtration.\newline
\no Set $\tau$ a stopping time.\newline

\no We assume that $\W$ is equipped with a probability measure $\nu$. For $p\geq 0$, we denote
$$L^p_a(\nu,H)=\left\{u\in L^p(\nu,H),\dot{u}\;is\;(\F_t)-adapted\right\}$$
\no and
$$\D=\left\{u\in L^0_a(\nu,H), \dot{u}\;is\;d\nu\times
  dt-a.s.\;bounded\right\}$$
\no For $t\in[0,1]$, we define
$$\pi_t:u\in L^0_a(\nu,H)\mapsto\int_0^.\dot{h}(s)1_{s\leq t}ds$$
\no Similarly, we define
\beaa\pi_\tau&:&u\in L^0_a(\nu,H)\mapsto\int_0^.\dot{h}(s)1_{s\leq \tau}ds\\
I-\pi_\tau&:&u\in L^0_a(\nu,H)\mapsto\int_0^.\dot{h}(s)1_{s> \tau}ds\eeaa
\no Notice that
\beaa\pi_\tau\D&\subset&\D\\
(I-\pi_\tau)\D&\subset&\D\eeaa
\no and define
$$\D_\tau=(I-\pi_\tau)\D$$
\no The filtration of a process m will be denoted
$\left(\F^m_t\right)$, the filtration of W will be simply denoted
$\left(\F_t\right)$. Except if stated otherwise, every filtration
considered is completed with respect to $\nu$.
\no If m is a martingale and v has a density whose stochastic
integral with respect to m is well defined we will denote
$$\delta_m v =\int_0^1\dot{v}(s)dm(s)$$
\no We also denote the Wick exponential as follow
$$\rho(\delta_m
v)=\exp\left(\int_0^1\dot{v}(s)dm(s)-\frac{1}{2}\int_0^1\left|\dot{v}(s)\right|^2d\langle
  m\rangle(s)\right)$$
\no and for $p\geq 0$ we denote
$$G_p(\nu,m)=\left\{u\in L^p_a(\nu,H), \EE_\nu\left[\rho(-\delta_m u)\right]=1\right\}$$
\no We assume there exists a family of adapted processes
$\left(W^u\right)_{u\in\D}$ and a $\nu$-Brownian motion $\beta$ which verify the following conditions:\newline
\no (i) $\beta$ is a $\nu$-Brownian motion whose canonical filtration
is identical to the canonical filtration of W\newline
\no (ii) $W^0=W$\newline
\no (iii) For every $u\in\D$, the law of $W^u$
under $\tilde{\nu}^u$ is the same as the law of W under $\nu$, where
$\tilde{\nu}^u$ is defined by
$\frac{d\tilde{\nu}^u}{d\nu}=\rho(-\delta_\beta u)$\newline
\no (iv) For every $u\in\D$,
$$\beta\circ W^u=\beta+u$$\newline
\no (v) For every $u,v\in\D$,
$$W^u\circ W^v=W^{v+u\circ W^v}\;\nu-a.s.$$
\no (vi) For every $u\in\D$
$$\left(W^u(s\wedge\tau),s\leq 1\right)=\left(W^{\pi_\tau u}(s\wedge\tau),s\leq 1\right)$$

\begin{remark}
\no Clearly $\D\subset L^\infty_a(\nu,H)$, so if $u\in\D$,
$\EE_\nu\left[\rho(-\delta_\beta u)\right]=1$ and $\tilde{\nu}^u$ which
was defined in condition (iii) is indeed a probability measure.
\end{remark}

\no Condition (iii) can be written as follow:

\begin{proposition}
\label{3gir}
\no Set $u\in \D$, for every bounded measurable function f, we have:
$$\EE_\nu\left[f\right]=\EE_\nu\left[f\circ W^u\rho(-\delta_\beta
  u)\right]$$
\end{proposition}

\no Next proposition ensures that the compositions written in (iv) and
(v) are well
defined.

\begin{proposition}
Set $u\in \D$, we have
$$W^u\nu\sim\nu$$
\end{proposition}

\nproof Set $f\in C_b(\W)$ bounded and measurable, we have, using
proposition \ref{3gir}
\beaa
\EE_{W^u\tilde{\nu}^u}\left[f\right]&=&\EE_{\tilde{\nu}^u}\left[f\circ W^u\right]\\
&=&\EE_{\nu}\left[f\circ
  W^u\rho\left(-\delta_\beta u\right)\right]\\
&=&\EE_{\nu}\left[f\right]\eeaa
\no so $W^u\tilde{\nu}^u=\nu$.\newline
\no Since $\tilde{\nu}\sim\nu$, we have $W^u\tilde{\nu}\sim W^u\nu$
which conclude the proof.\nqed

\begin{definition}
\label{3DD}
 Set $\DD$ a subset of $G_0(\nu,\beta)$ such that the map
$u\in\D\mapsto W^u$ can be extended to $\DD$ while verifying the following conditions.\newline
\no (i) $\D\subset\DD\subset G_2(\nu,\beta)$\newline
\no (ii) For any $u\in\DD$, $W^u$ is adapted.\newline
\no (iii)  For every $u\in\DD$, the law of $W^u$
under $\tilde{\nu}^u$ is the same as the law of W under $\nu$, where
$\tilde{\nu}^u$ is defined by
$\frac{d\tilde{\nu}}{d\nu}=\rho(-\delta_\beta\nu)$\newline
\no (iv) For every $u\in\DD$,
$$\beta\circ W^u=\beta+u$$\newline
\no (v) For every $u,v\in\DD$ such that $v+u\circ W^v\in\DD$
$$W^u\circ W^v=W^{v+u\circ W^v}\;\nu-a.s.$$\newline
\no (vi) There exists $\DDD$ such that $D''\subset \DDD\subset
L^0_a(\nu,H)$, $\DD=\DDD\cap G_2(\nu,\beta)$ and for every $u\in\DD$ such that the equation $u+v\circ W^u$ has
a solution in $G_0(\nu,\beta)$, this equation has a solution in
$\DDD$.\newline
\no (vii) For every $u\in\DD$ such that $\pi_\tau u\in\DD$,
$$\left(W^u(s\wedge\tau),s\leq 1\right)=\left(W^{\pi_\tau u}(s\wedge\tau),s\leq 1\right)$$\newline
\end{definition}

\no\begin{remark}
$\D$ verify the set of condition above.
\end{remark}\newline

\begin{proposition}
Set $u\in\DD$. For every bounded measurable function f, we have
$$\EE_\nu\left[f\right]=\EE_\nu\left[f\circ W^u\rho(-\delta_\beta u)\right]$$
\no Furthermore,
$$W^u\nu\sim\nu$$
\end{proposition}

\nproof The first assertion is condition (iii). The proof of the second assertion is the same as the case $u\in\D$.\nqed

\begin{definition}
We define $\DD_\tau$ as
$$\DD_\tau=\DD\cap (I-\pi_\tau)L^0_a(\nu,H)$$
\end{definition}

\section{\bf{Conditional expectation results}}

\no We need the abstract Bayes formula for a stopping time:

\begin{lemma}
Set $\theta$ a probability measure on $\left(\W,\F\right)$ such that $\theta\ll\nu$ . Denote
$$L=\frac{d\theta}{d\nu}$$
\no For every measurable $f:\W\rightarrow \RR$ we have
$$\EE_\theta\left[f|\F_\tau\right]=\frac{\EE_\nu\left[fL|\F_\tau\right]}{\EE_\nu\left[L|\F_\tau\right]}$$
\end{lemma}

\nproof We can assume f is positive. Denote for $s\in [0,1]$
$$L(s)=\EE_\nu\left[L|\F_s\right]$$
\no The martingale stopping theorem gives
$$L(\tau)=\EE_\nu\left[L|\F_\tau\right]$$
\no Set $A\in\F_\tau$, we need to prove that
$$\EE_\nu\left[1_AL(\tau)\EE_\theta\left[f|\F_\tau\right]\right]=\EE_\nu\left[1_A\EE_\nu\left[fL|\F_\tau\right]\right]$$
\no We have
\beaa \EE_\nu\left[1_AL(\tau)\EE_\theta\left[f|\F_\tau\right]\right]&=&\EE_{\nu|_{\F_\tau}}\left[1_AL(\tau)\EE_\theta\left[f|\F_\tau\right]\right]\\
&=&\EE_{\theta|_{\F_\tau}}\left[1_A\EE_\theta\left[f|\F_\tau\right]\right]\\
&=&\EE_{\theta}\left[1_A\EE_\theta\left[f|\F_\tau\right]\right]\\
&=&\EE_{\theta}\left[1_Af\right]\\
&=&\EE_{\nu}\left[1_AfL\right]\\
&=&\EE_{\nu}\left[1_A\EE_\nu\left[fL|\F_\tau\right]\right]\eeaa\nqed

\begin{proposition}
Set $u\in \DD$ and $f\in
L^0(\nu)$ an $\F_\tau$-measurable function. Then $\nu$-a.s.
$$f\circ  W^u=f\circ  W^{\pi_\tau u}$$
\no Consequently, if $u\in \DD_\tau$, we have $\nu$-a.s.
$$f\circ  W^u=f$$
\end{proposition}
\nproof For $s\in [0,1]$, we have
\beaa\beta(s)\circ W^u&=&\beta(s)+u(s)\\
&=&\beta(s)+\pi_s u(s)\\
&=&\beta(s)\circ W^{\pi_s u}\eeaa
\no Consequently, for $h\in H$
$$\rho(\delta_\beta \pi_s h)\circ W^u= \rho(\delta_\beta \pi_s h)\circ W^{\pi_s u}$$
\no and
$$\rho(\delta_\beta \pi_\tau h)\circ W^u= \rho(\delta_\beta \pi_\tau h)\circ W^{\pi_\tau u}$$
\no Denote $$L^2(\nu,\F_\tau)=\{f\in L^2(\nu),f\;is\;\F_\tau-measurable\}$$
\no $\left(\rho(\delta_\beta \pi_s h),s\in [0
,1]\right)$ being a closed martingale, we have
$$\EE_\nu\left[\rho(\delta_\beta  h)|\F_\tau\right]=\rho(\delta_\beta \pi_\tau h)$$
\no Since $\beta$ and W have the same filtration, the vector space generated by $\{\rho(\delta_\beta h),h\in H\}$ is dense in $L^2(\nu)$. $g\in L^2(\nu)\mapsto \EE_\nu\left[g|\F_\tau\right]$ being a continuous surjection from $L^2(\nu)$ to $L^2(\nu,\F_\tau)$, the vector space generated by
$\{\rho(\delta_\beta \pi_\tau h),h\in H\}$ is dense in $L^2(\nu,\F_\tau)$. We denote E this vector space.\newline
\no Assume that f is bounded, there exists $(f_n)\in E^\NN$ which converges to f
$\nu$-a.s. Since $W^u\nu\ll\nu$ and
$W^{\pi_\tau u}\nu\ll\nu$, $\left(f_n\circ  W^u\right)$
converges $\nu$-a.s. to $f\circ  W^u$ and $\left(f_n\circ W^{\pi_\tau u}\right)$
converges $\nu$-a.s. to $f\circ  W^{\pi_\tau u}$, which
ensures the result in this case\newline
\no Finally, if f is only supposed to be $\F_\tau$-measurable, there
exists a sequence of bounded $\F_\tau$-measurable functions which
converges to f and we proceed as above.\nqed

\begin{proposition}
Set L a density on $(\W,\nu,\F)$ such that $L>0$
$\nu$-a.s. Denote,
$$M(s)=\EE_{\nu}\left[L\left|\F_s\right.\right]$$
\no and set $v\in
L^0_a(\nu,H)$ such that $$M(s)=\rho(-\delta_\beta\pi_sv)$$
\no Then the two following propositions are equivalent:\newline
\no (i) $\pi_\tau v=0$ $\nu$-a.s.\newline
\no (ii) $\EE_{\nu}\left[L\left|\F_\tau\right.\right]=1$ $\nu$-a.s.
\end{proposition}
\nproof The direct implication is trivial. Conversely, M is a martingale with unit expectation and since
$$M(s)=1+\int_0^sM(r)\dot{v}(r)d\beta(r)$$
\no We have
$$\left\langle M-1\right\rangle = \int_0^.\left(M(r)\dot{v}(r)\right)^2dr$$
\no Proposition(ii) gives $\left(M(s\wedge\tau)-1,s\leq 1\right)=0$ $\nu$-a.s., so $\left(\left\langle M-1\right\rangle(s\wedge\tau),s\leq 1\right)=0$ $\nu$-a.s., $L>0$ so $\nu$-a.s. $M(s)>0$ for any $s\in[0,1]$ and we have proposition (i).\nqed

\begin{lemma}
Set $u\in \DD_\tau$ and denote $L=\frac{d W^u\nu}{d\nu}$. We have
$$\EE_{\nu}\left[L\left|\F_\tau\right.\right]=1$$
\no Consequently, for any $f\in L^1\left( W^u\nu\right)$, we have
$$\EE_{ W^u\nu}\left[f\left|\F_\tau\right.\right]=\EE_{\nu}\left[f\circ W^u\left|\F_\tau\right.\right]=\EE_{\nu}\left[fL\left|\F_\tau\right.\right]$$
\end{lemma}
\nproof Set $B\in\F_\tau$, we have
\beaa \EE_{\nu}\left[1_BL\right]=\EE_{\nu}\left[1_B\circ W^u\right]=\EE_{\nu}\left[1_B\circ W^{\pi_\tau u}\right]=\EE_{\nu}\left[1_B\right]\eeaa
\no Then the second assertion is a direct consequence of Bayes
formula.\nqed

\section{\bf{Invertibility results}}

\begin{definition}
A measurable map $U:\W\rightarrow \W$ is said to be
$\nu$-a.s. left-invertible if and only if $U\nu\ll\nu$ and there exists a measurable map
$V:\W\rightarrow \W$ such that $V\circ U=I_\W$ $\nu$-a.s.\newline
\no A measurable map $U:\W\rightarrow \W$ is said to be
$\nu$-a.s. right-invertible if and only if there exists a measurable map
$V:\W\rightarrow \W$ such that $V\nu\ll\nu$ and  $U\circ V=I_\W$ $\nu$-a.s.
\end{definition}

\begin{proposition}
\label{3lrinv}
Set $U,V:\W\rightarrow\W$ measurable maps such that $V\circ U=I_\W$
$\nu$-a.s. and $V\nu\ll\nu$ Then  $U\circ V=I_\W$ $U\nu$-a.s., so if $U\nu\sim \nu$, we
also have $U\circ V=I_\W$ $\nu$-a.s. In that case, we will say that
$U$ is $\nu$-a.s. invertible and we also have $V\nu\sim\nu$.
\end{proposition}

\nproof See \cite{a2}.\nqed

\begin{proposition}
\no Set
$u\in\DD_\tau\cap L^2_a(\nu,H)$. If $W^u$ is $\nu$-a.s. left-invertible. Then there exists
$v\in\DD_\tau$ such that $\nu$-a.s.
$$W^v\circ W^u=W^u\circ W^v=I_\W$$
\no and
\beaa \frac{dW^u\nu}{d\nu}&=&\rho(-\delta_\beta v)\\
\frac{dW^v\nu}{d\nu}&=&\rho(-\delta_\beta u)\eeaa
\end{proposition}
\nproof Everything is already known from \cite{a2} except the fact that $\pi_\tau v =0$. This arises from the relation
$$\dot{v}(s)=-\dot{u}(s)\circ  W^v$$\nqed

\no Now we recall two very useful lemmas, see \cite{a2} for the proof

\begin{lemma}
\label{3li1}
 Set $u\in\DD\cap L^2_a(\nu,H)$ and denote
 $L=\frac{d W^u\nu}{d\nu}$, we have $\nu$-a.s.
$$L\circ  W^u\EE_{\nu}\left[\rho(-\delta_\beta u)\left|\F^{ W^u}_1\right.\right]= 1$$
\end{lemma}

\begin{theorem}
\label{3pi2}
Set $u\in\DD\cap L^2_a(\nu,H)$ and denote
 $L=\frac{d W^u\nu}{d\nu}$. Then $ W^u$ is $\nu$-a.s. left-invertible if and only if
$$\EE_{\nu}\left[L\log L\right]=\frac{1}{2}\EE_{\nu}\left[|u|^2_H\right]$$
\no Moreover, if $ W^u$ is $\nu$-a.s. left-invertible, we have $\nu$-a.s.
$$L\circ  W^u\rho\left(-\delta_\beta u\right)=1$$
\end{theorem}

\no Now we give the results relative to the invertibility of $ W^u$ when $\pi_\tau u=0$.

\begin{proposition}
\no Set $u\in\DD_\tau\cap L^2_a(\nu,H)$ and denote $L=\frac{d W^u\nu}{d\nu}$. We have $\nu$-a.s.:
$$\EE_{\nu}\left[L\log L\left|\F_\tau\right.\right]\leq\frac{1}{2}\EE_{\nu}\left[|u|^2_H\left|\F_\tau\right.\right]$$
\end{proposition}

\nproof We have $\left( W^u(s\wedge\tau)\right)_{s\leq 1}=\left( W(s\wedge\tau)\right)_{s\leq 1}$ hence
$$\F_\tau=\F^{ W^u}_\tau\subset \F^{ W^u}_1$$
\no Consequently, using lemma \ref{3li1} and Jensen inequality, we have:
\beaa \EE_{\nu}\left[L\log L\left|\F_\tau\right.\right]&=&\EE_{\nu}\left[\log L\circ W^u\left|\F_\tau\right.\right]\\
&\leq&-\EE_{\nu}\left[\left.\log\EE_{\nu}\left[\rho(-\delta_\beta u)\left|\F^{ W^u}_1\right.\right]\right|\F_\tau\right]\\
&\leq&-\EE_{\nu}\left[\left.\EE_{\nu}\left[\log\rho(-\delta_\beta u)\left|\F^{ W^u}_1\right.\right]\right|\F_\tau\right]\\
&\leq&-\EE_{\nu}\left[\log\rho(-\delta_\beta u)\left|\F_\tau\right.\right]\\
&\leq&\frac{1}{2}\EE_{\nu}\left[|u|_H^2\left|\F_\tau\right.\right]\eeaa\nqed

\begin{theorem}
\label{3inv}
Set $u\in \DD_\tau\cap L^2_a(\nu,H)$ and denote
 $L=\frac{d W^u\nu}{d\nu}$. Then $ W^u$ is $\nu$-a.s. left-invertible,
if and only if $\nu$-a.s.
$$\EE_{\nu}\left[L\log L\left|\F_\tau\right.\right]=\frac{1}{2}\EE_{\nu}\left[|u|^2_H\left|\F_\tau\right.\right]$$
\end{theorem}

\nproof Assume that the equality holds, taking the expectation we have
$$\EE_{\nu}\left[L\log L\right]=\frac{1}{2}\EE_{\nu}\left[|u|^2_H\right]$$
\no so according to theorem \ref{3pi2} $W^u$ is $\nu$-a.s. left invertible.\newline
\no Conversely, using again theorem  \ref{3pi2}, we have
$$L\circ  W^u\rho\left(-\delta_\beta u\right)=1$$
\no So, since $\pi_\tau u=0$,
\beaa \EE_{\nu}\left[L\log L\left|\F_\tau\right.\right]&=&\EE_{\nu}\left[\log L\circ W^u\left|\F_\tau\right.\right]\\
&=&\EE_{\nu}\left[-\log \rho(-\delta_\beta u)\left|\F_\tau\right.\right]\\
&=&\frac{1}{2}\EE_{\nu}\left[|u|^2_H\left|\F_\tau\right.\right]\eeaa\nqed

\begin{definition}
We denote
\beaa \D^i&=&\left\{u\in \D,W^u\;is\;\nu-a.s.\;invertible\right\}\\
\D_\tau^i&=&\D_\tau\cap\D^i\eeaa
\end{definition}

\section{\bf{Approximation of absolutely continuous measures}}

\begin{theorem}
\label{3repr}
If $\theta\sim \nu$ is such that there exists $p>1$ such that
$$\frac{d\theta}{d\nu}\in L^p(\nu)$$
\no and $\nu$-a.s.
$$\EE_{\nu}\left[\left.\frac{d\theta}{d\nu}\right|\F_\tau\right]=1$$
\no There exists $(u_n)\in\left(\D_\tau^i\right)^\NN$ such that,
\beaa \frac{d W^{u_n}\nu}{d\nu}\rightarrow \frac{d\theta}{d\nu}\;\; in\;L^p(\nu)\eeaa
\end{theorem}

\nproof \no Eventually sequentializing afterward, we have to prove that for any
$\epsilon>0$, there exists $u\in \D_\tau^i$ such that
\beaa \left|\frac{d W^u\nu}{d\nu}-\frac{d\theta}{d\nu}\right|_{L^p(\nu)}&\leq&\epsilon\eeaa
\no The proof is divided in six steps.\newline
\no Step 1 : We approximate $\frac{d\theta}{d\nu}$ with a density that is both lower and
upper bounded.\newline
\no Denote
$$L(s)=\EE_{\nu}\left[\left. \frac{d\theta}{d\nu}\right|\F_s\right]$$
\no  and for $n\in\NN$,
$$T_n=\inf\left\{s\in[0,1], L(s)\geq n\right\}$$
\no $L(1)=\frac{d\theta}{d\nu}$ and L being a closed martingale, $L^{T_n}(.)$ is still a closed martingale which converges in $L^1$ to $L(T_n)=\EE_{\nu}\left[L(1)\left|\F_{T_n}\right.\right]$, so
\beaa \EE_{\nu}\left[L(T_n)|\F_\tau\right]&=&\EE_{\nu}\left[\left.L^{T_n}(1)\right|\F_\tau\right]\\
&=&L^{T_n}(\tau)\\
&=&L(\tau\wedge T_n)\\
&=&\EE_{\nu}\left[L(1)\left|\F_{\tau\wedge T_n}\right.\right]\\
&=&\EE_{\nu}\left[\EE_{\nu}\left[L(1)\left|\F_{\tau}\right.\right]\left|\F_{\tau\wedge T_n}\right.\right]\\
&=&1\eeaa
\no Furthermore, since L is a closed martingale, $\left(\EE_{\nu}\left[L(1)\left|\F_{T_n}\right.\right]\right)_{n\in\NN}$ is also a closed martingale so is uniformly integrable. $\left(L(T_n)\right)$ converges to $L_1$ $\nu$-a.s. and Jensen inequality gives
\beaa 0\leq L(T_n)^p=\EE_{\nu}\left[L(1)\left|\F_{T_n}\right.\right]^p\leq\EE_{\nu}\left[L(1)^p\left|\F_{T_n}\right.\right]\eeaa
\no So $\left(L(T_n)^p\right)$ is uniformly integrable and $\left(L(T_n)\right)$ converges in $L^p(\nu)$ to $L(1)$, so there exists $n_0\in\NN$ such that
$$\left|L(T_{n_0})-L(1)\right|_{L^r(\nu)}\leq\epsilon$$

\no $\left(\frac{L(T_{n_0})+a}{1+a}\right)$ converges
$\nu$-a.s. to $L(T_{n_0})$ when a converges to 0. Set $a\in [0,1]$, we have
\beaa 0\leq\frac{L(T_{n_0})+a}{1+a}\leq L(T_{n_0})+1\eeaa
\no $L(T_{n_0})+1\in L^p(\nu)$ so according to the Lebesgue theorem, $\left(\frac{L(T_{n_0})+a}{1+a}\right)$ converges
to $L(T_{n_0})$ in $L^p(\nu)$ and
there exists $a\in [0,1]$ such that
\beaa \left|\frac{L(T_{n_0})+a}{1+a}-L(T_{n_0})\right|_{L^p(\nu)}\leq\epsilon\eeaa

\no $\frac{L(T_{n_0})+a}{1+a}$ is both lower-bounded and upper-bounded in
$L^\infty(\nu)$,
denote these bounds respectively d and D.\newline
\no Denote
$$M(s)=\EE_{\nu}\left[\left. \frac{L(T_{n_0})+a}{1+a}\right|\F_s\right]$$
\no We can write
$$M=\exp\left(\int_0^.\dot{\alpha}(r)d\beta(r)-\frac{1}{2}\int_0^.\left|\dot{\alpha}(r)\right|^2dr\right)$$
\no with $\alpha\in
(I-\pi_\tau) L^0_a(\nu,H)$ since
$$\EE_{\nu}\left[M(1)|\F_\tau\right]=\frac{\EE_{\nu}\left[\left.L_{T_{n_0}}\right|\F_\tau\right]+a}{1+a}=1$$
\no Step 2 : we prove that $\alpha\in (I-\pi_\tau) L^2_a(\nu,H)$\newline
\no Set
$$S_n=\inf\left\{s\in [0,1], \int_0^s\left|\dot{\alpha}(r)\right|^2dr>n\right\}$$
\no $(S_n)$ is a sequence of stopping times which increases
stationarily toward 1. We have, using
$M=1+\int_0^.\dot{\alpha}(r)M(r)d\beta(r)$
\beaa \EE_{\nu}\left[\left(M(s\wedge
      S_n)-1\right)^2\right]&=&\EE_{\nu}\left[\int_0^{s\wedge S_n}\left|\dot{\alpha}(r)\right|^2M(r)^2dr\right]\\
&\geq& d^2\EE_{\nu}\left[\int_0^{s\wedge S_n}\left|\dot{\alpha}(r)\right|^2dr\right]\eeaa
\no so
$$\EE_{\nu}\left[\int_0^{s\wedge S_n}\left|\dot{\alpha}(r)\right|^2dr\right]\leq \frac{1}{d^2}\EE_{\nu}\left[\left(M(s\wedge
      S_n)-1\right)^2\right]\leq \frac{2\left(D^2+1\right)}{d^2}$$
\no hence passing to the limit
$$\EE_{\nu}\left[\int_0^1\left|\dot{\alpha}(r)\right|^2dr\right]\leq\infty$$
\no Step 3 : We approximate $\alpha$ with an element of $(I-\pi_\tau) L^\infty_a(\nu,H)$.\newline
\no Define
$$\alpha^n(s,w)\in [0,1]\times \W\mapsto\int_0^s
\dot{\alpha}(r,w)1_{[0,S_n]}(r,w)dr$$

\no and
\beaa M^n(s)&=&\exp\left(\int_0^s\dot{\alpha^n}(r)d\beta(r)-\frac{1}{2}\int_0^s\left|\dot{\alpha^n}(r)\right|^2dr\right)\eeaa
\no $\alpha^n\in (I-\pi_\tau) L^\infty(\nu,H)$ and $M^n(1)=\EE_{\nu}\left[M(1)\left|\F_{S_n}\right.\right]$ so $(M^n(1))_{n\in\NN}$ is a closed martingale since M is one, hence it converges $\nu$-a.s. to $M(1)$ and it is uniformly integrable. Jensen inequality gives
\beaa 0\leq\left|M^n(1)\right|^p\leq\EE_{\nu}\left[M(1)\left|\F_{S_n}\right.\right]^p\leq\EE_{\nu}\left[M(1)^p\left|\F_{S_n}\right.\right]\eeaa
\no So $\left(\left|M^n(1)\right|^r\right)$ is uniformly integrable and $(M^n(1))$ converges to $M(1)$ in $L^p(\nu)$. Consequently, there exists $n\in\NN$ such that

\beaa \left|M^{n}(1)-M(1)\right|_{L^p(\nu)}&\leq&\epsilon\eeaa
\no Step 4 : we approximate $\alpha^n$ with an element of $\D_\tau$

\no Define
$$\xi^{n,m}:(s,w)\in [0,1]\times\W\mapsto\int_0^s\max \left(\min
  \left(\dot{\alpha}^n(r,w),m\right),-m\right)dr$$
\no and
$$
M^{n,m}(s)=\exp\left(\int_0^s\dot{\xi^{n,m}}(r)d\beta(s)-\frac{1}{2}\int_0^s\left|\dot{\xi^{n,m}}(r)\right|^2dr\right)$$
\no $\xi^{m,n}\in\D_\tau$ and $\left(M^{n,m}(1)\right)$ and converges to $M^n(1)$ in probability. To prove
that $\left((M^{n,m}(1)^p\right)$ is uniformly integrable, it is
sufficient to prove that $\left(M^{n,m}(1)\right)$  is bounded in every
$L^q(\nu)$ with $q>1$. Set $q>1$
\beaa
\EE_\nu\left[\left|M^{n,m}_1\right|^{q}\right]&=&\EE_\nu\left[\exp\left(q\int_0^1\dot{\xi^{n,m}}(s)d\beta(s)-\frac{q}{2}\int_0^1\left|\dot{\xi^{n,m}}(s)\right|^2ds\right)\right]\\
&=&\EE_\nu\left[\exp\left(q\int_0^1\dot{\xi^{n,m}}(s)d\beta(s)-\frac{q^2}{2}\int_0^1\left|\dot{\xi^{n,m}}(s)\right|^2ds\right)\exp\left(\frac{q^2-q}{2}\int_0^1\left|\dot{\xi^{n,m}}(s)\right|^2ds\right)\right]\\
&\leq&\EE_\nu\left[\exp\left(\int_0^1q\dot{\xi^{n,m}}(s)d\beta(s)-\frac{1}{2}\int_0^1\left|p\dot{\xi^{n,m}}(s)\right|^2ds\right)\exp\left(\frac{q^2-q}{2}n\right)\right]\\
&\leq&\exp\left(\frac{q^2-q}{2}n\right)\eeaa
\no so $\left(M^{n,m}(1),m\in\NN\right)$ converges to $M_1^n$ in $L^p(\nu)$ and there exists some $m>0$ such that
\beaa \left|M^{n,m}(1)-M^n(1)\right|_{L^p(\nu)}&\leq&\epsilon\eeaa

\no Step 5 : We approximate $\xi^{n,m}$ with a retarded shift
$\gamma^\eta$, so that $ W^{\gamma^\eta}$ is
$\nu$-a.s. invertible.\newline
\no For a given $\eta>0$, set
$$\gamma^{\eta}(s,w)\in[0,1]\times W\mapsto\int_0^s
\dot{\xi^{n,m}}(r-\eta)1_{r\geq \eta}ds$$
\no and
\beaa N^{\eta}(s)&=&\exp\left(\int_0^s\dot{\gamma^\eta}(r)d\beta(r)-\frac{1}{2}\int_0^s\left|\dot{\gamma^\eta}(r)\right|^2dr\right)\eeaa
\no Clearly $\gamma^\eta\in\D_\tau$ and $\gamma^{\eta}\rightarrow \xi^{n,m}$ in
$L^2(\nu,H)$ when
$\eta\rightarrow 0$, which ensures that
$\left(N^\eta(1),\eta>0\right)$ converges to $M^{n,m}(1)$ in
probability.\newline
\no As in step 4, $\left(N^\eta(1),\eta>0\right)$ is bounded in every
$L^q(\nu)$ and so $\left(N^\eta(1)^p,\eta>0\right)$ is uniformly integrable
and $\left(N^\eta(1),\eta>0\right)$  converges to $M^{n,m}(1)$
in $L^p(\nu)$. There exists $\eta>0$ such that
\beaa \left|N^{\eta}(1)-M^{n,m}(1)\right|_{L^p(\nu)}&\leq&\epsilon\eeaa
\no Using triangular inequality, we have
\beaa \left|\frac{d\theta}{d\nu}-N^\eta(1)\right|_{L^p(\nu)}&\leq&\left|\frac{d\theta}{d\nu}-L(T_{n_0})\right|_{L^p(\nu)}+\left|L(T_{n_0})-\frac{L(T_{n_0})+a}{1+a}\right|_{L^p(\nu)}\\
&&+\left|\frac{L_(T_{n_0})+a}{1+a}-M^n(1)\right|_{L^p(\nu)}\\
&&+\left|M^{n}(1)-M^{n,m}(1)\right|_{L^p(\nu)}\\
&&+\left|M^{n,m}(1)-N^\eta(1)\right|_{L^p(\nu)}\\
&\leq&5\epsilon\eeaa

\no Step 6 : We prove that $ W^{-\gamma^{\eta}}$ is
$\nu$-a.s. left-invertible and is the solution to our problem.\newline
\no Set $A\subset\W$ such that $\nu(A)=1$ and for every $w\in A$, $\beta\circ W^{-\gamma^\eta}(w)=\beta(w)-\gamma^\eta(w)$ and set $w_1,w_2\in A$ such that
$ W^{-\gamma^\eta}(w_1)= W^{-\gamma^\eta}(w_2)$. We have
\beaa \beta\circ W^{-\gamma^\eta}(w_1)&=&\beta\circ W^{-\gamma^\eta}(w_2)\\
\beta(w_1)-\int_0^.\dot{\gamma^\eta}(s,w_1)ds&=&\beta(w_2)-\int_0^.\dot{\gamma^\eta}(s,w_2)ds\eeaa

\no For any $s\in[0,\eta]$, $\beta(s,w_1)=\beta(s,w_2)$,
$\gamma^\eta$ being adapted to filtration $(\F^\beta_{s-\eta})$,
it implies that for $s\in[0,2\eta]$
$$\int_0^s\dot{\gamma^\eta}(r,w_1)ds=\int_0^s\dot{\gamma^\eta}(r,w_2)ds$$
\no and
$$\beta(s,w_1)=\beta(s,w_2)$$
\no An easy iteration shows that $\beta(w_1)=\beta(w_2)$.\newline

\no Since $\beta$ and W have the same filtrations and $\beta$ is $\mu$-a.s. path-continuous, we can write $W(t)=\phi_t(\beta(s),s\in [0,t]\cap\QQ)$ $\nu$-a.s. for every $t\in [0,1]$, with $\phi_t$ a measurable function from $\RR^\QQ$ to $\RR$, see \cite{nev}. Consequently, we can write $\left(W(t),t\in[0,1]\cap\QQ\right)=\phi\left(\beta(t), t\in[0,1]\cap\QQ\right)$ $\nu$-a.s., with $\phi$ a measurable function from $\RR^\QQ$ to $\RR^\QQ$. Denote
\beaa A'=A\cap\left\{w\in\W, \left(W(t,w),t\in[0,1]\cap\QQ\right)=\phi\left(\beta(t,w), t\in[0,1]\cap\QQ\right)\right\}\eeaa
\no $\nu(A')=1$. Set $w_1,w_2\in A'$ such that $ W^{-\gamma^\eta}(w_1)= W^{-\gamma^\eta}(w_2)$. We have $\beta(w_1)=\beta(w_2)$ so
\beaa\left(W(t,w_1),t\in[0,1]\cap\QQ\right)&=&\left(W(t,w_2),t\in[0,1]\cap\QQ\right)\\
\left(w_1(t),t\in[0,1]\cap\QQ\right)&=&\left(w_2(t),t\in[0,1]\cap\QQ\right)\eeaa

\no $w_1$ and $w_2$ are continuous and coincide on $[0,1]\cap\QQ$ so they are equal.

\no $ W^{-\gamma^{\eta}}$ is $\nu$-a.s. injective and so
$\nu$-a.s. left-invertible, its inverse is of
the form $ W^{v^{\eta}}$, with $v^\eta\in\D_\tau$ and we have
$$\frac{d  W^{v^\eta}\nu}{d\nu}=L^{\eta,n}_1$$
\no So $ W^{v^\eta}\nu\sim \nu$ and
$$ W^{v^\eta}\circ
 W^{-\gamma^\eta}=  W^{-\gamma^\eta}\circ  W^{v^\eta}\;\;
\nu-a.s.$$
\nqed

\begin{corollary}
\label{3repr2}
If $\theta\sim \nu$ is such that there exists
$q,p>1$ such that $p^{-1}+q^{-1}=1$ and
\beaa \frac{d\theta}{d\nu}&\in& L^p(\nu)\\
\log\frac{d\theta}{d\nu}&\in& L^q(\nu)\eeaa
\no and $\nu$-a.s.
$$\EE_{\nu}\left[\left.\frac{d\theta}{d\nu}\right|\F_\tau\right]=1$$
\no there exists $(u_n)\in \left(\D_\tau^i\right)^\NN$ such that
\beaa
\EE_{\nu}\left[\left.\frac{d W^{u_n}\nu}{d\nu}\log
    \frac{d W^{u_n}\nu}{d\nu}\right|\F_\tau\right]&\rightarrow&
\EE_{\nu}\left[\left.\frac{d\theta}{d\nu}\log
    \frac{d\theta}{d\nu}\right|\F_\tau\right]\;\;\nu-a.s.\\
\EE_{\nu}\left[\left.\frac{d W^{u_n}\nu}{d\nu}\log
    \frac{d\theta}{d\nu}\right|\F_\tau\right]&\rightarrow&
\EE_{\nu}\left[\left.\frac{d\theta}{d\nu}\log
    \frac{d\theta}{d\nu}\right|\F_\tau\right]\;\;\nu-a.s.\eeaa
\end{corollary}

\nproof From theorem \ref{3repr}, there exists $(u_n)\in
\left(\D_\tau^i\right)^\NN$ such that for
every n,
\beaa
&&\frac{d W^{u_n}\nu}{d\nu}\rightarrow
\frac{d\theta}{d\nu}\;\; in\;L^r(\nu)\eeaa
\no This implies
$$\frac{d W^{u_n}\nu}{d\nu}\log \frac{d W^{u_n}\nu}{d\nu} \rightarrow
\frac{d\theta}{d\nu}\log \frac{d\theta}{d\nu}\;\;
in\;L^1(\nu)$$
\no Holder inequality gives
\beaa \left|\frac{d W^{u_n}\nu}{d\nu}\log \frac{d\theta}{d\nu}-\frac{d\theta}{d\nu}\log \frac{d\theta}{d\nu}\right|_{L^1(\nu)}&\leq&\left|\frac{d W^{u_n}\nu}{d\nu}-
\frac{d\theta}{d\nu}\right|_{L^p(\nu)}\left|\log
\frac{d\theta}{d\nu}\right|_{L^q(\nu)}\\
&\rightarrow&0\eeaa
\no The corresponding conditional expectations converges similarly in
$L^1(\nu)$ since $\EE_{\nu}\left[.|\F_\tau\right]$ is a
bounded operator with norm 1 in $L^1(\nu)$. Finally we can
extract a subsequence of $(u_n)$ to get the two desired almost sure convergences.\nqed

\section{\bf{Variational problem}}
\no As stated in the beginning, we aim to provide a variational
representation of $-\log\EE_{\nu}\left[\left.e^{-f}\right|\F_\tau\right]$.

\begin{definition}
We denote $\mathcal{P}_\tau$ the set of probability measures $\theta$ on
$\left(\W,\F\right)$ such that \beaa
\theta&\sim&\nu\\
\EE_{\nu}\left[\left.\frac{d\theta}{d\nu}\right|\F_\tau\right]&=&1\eeaa
\end{definition}

\begin{theorem}
Set $f:\W\rightarrow \RR$ a measurable function verifying
$$\EE_{\nu}\left[|f|(1+e^{-f})\right]<\infty$$
\no Then
$$-\log\EE_{\nu}\left[\left.e^{-f}\right|\F_\tau\right]=\inf_{\theta\in\mathcal{P}_\tau}\EE_\theta\left[\left.f+\log\frac{d\theta}{d\nu}\right|\F_\tau\right]\;\;\nu-a.s.$$
\no and the unique infimum is attained at the measure
$$d\theta_0=\frac{e^{-f}}{\EE_{\nu}\left[\left.e^{-f}\right|\F_\tau\right]}d\nu$$
\end{theorem}

\nproof Set $\theta\in\mathcal{P}_\tau$, denote
$$L(s)=\left.\frac{d\theta}{d\nu}\right|_{\F_s}$$
\no $L(\tau)=1$ $\nu$-a.s. since $\theta\in\mathcal{P}_\tau$ so using the Bayes formula:
\beaa\log\EE_\nu\left[e^{-f}|\F_\tau\right]&=&\log\EE_\nu\left[\left.e^{-f}\frac{L(1)}{L(1)}\right|\F_\tau\right]\\
&=&\log\EE_\theta\left[\left.e^{-f}\frac{L(\tau)}{L(1)}\right|\F_\tau\right]\\
&=&\log\EE_\theta\left[\left.\frac{e^{-f}}{L(1)}\right|\F_\tau\right]\eeaa
\no Jensen inequality gives
\beaa -\log\EE_\theta\left[\left.\frac{e^{-f}}{L(1)}\right|\F_\tau\right]&\leq&\EE_\theta\left[\left.-\log\frac{e^{-f}}{L(1)}\right|\F_\tau\right]\\
&\leq&\EE_\theta\left[f|\F_\tau\right]+\EE_\theta\left[\log L(1)|\F_\tau\right]\\
&\leq&\EE_\theta\left[f|\F_\tau\right]+\EE_\theta\left[\log L(1)|\F_\tau\right]\eeaa
\no A straightforward calculation gives
$$\EE_{\theta_0}\left[\left.f+\log\frac{d\theta_0}{d\nu}\right|\F_\tau\right]=-\log\EE_\nu\left[e^{-f}|\F_\tau\right]$$
\no and the reverse inequality.\nqed

\begin{proposition}
\label{3ineqvar}
\no Set $f:\W\rightarrow\RR$ a measurable function
verifying $\EE_{\nu}\left[|f|(1+e^{-f})\right]<\infty$, then
$$-\log\EE_{\nu}\left[\left.e^{-f}\right|\F_\tau\right]\leq\inf_{u\in \DD_\tau\cap L^2_a(\nu,H)}\EE_{\nu}\left[\left.f\circ W^u+\frac{1}{2}|u|_H^2\right|\F_\tau\right]\;\;\nu-a.s.$$
\end{proposition}

\nproof Denote
 $\mathcal{P}_\tau'$ the set of the elements S of $\mathcal{P}_\tau$ such that
 there exists some $u\in \DD_\tau$ which verifies
 $S= W^u\nu$.\newline
\no Set $\theta\in\mathcal{P}'_\tau$ and denote
$L=\frac{d\theta}{d\nu}$. Since $\EE_{\nu}\left[L|\F_\tau\right]=1$, we
have using Bayes formula
\beaa
\EE_\theta\left[f|\F_\tau\right]&=&\EE_{\nu}\left[fL|\F_\tau\right]\\
&=&\EE_{\nu}\left[\left.f\circ W^u\right|\F_\tau\right]\\
\EE_\theta\left[\log L|\F_\tau\right]&=&\EE_\nu\left[L\log
  L|\F_\tau\right]\\
&\leq&\frac{1}{2}\EE_{\nu}\left[\left.|u|^2_H\right|\F_\tau\right]\eeaa
\no So since $\mathcal{P}_\tau'\subset\mathcal{P}_\tau$, we have
\beaa
-\log\EE_{\nu}\left[\left.e^{-f}\right|\F_\tau\right]&=&\inf_{\theta\in\mathcal{P}_\tau}\left(\EE_\theta[f|\F_\tau]+\EE_\theta\left[\left.\log\frac{d\theta}{d\nu}\right|\F_\tau\right]\right)\\
&\leq&\inf_{\theta\in\mathcal{P}_\tau'}\left(\EE_\gamma[f|\F_\tau]+\EE_\theta\left[\left.\log\frac{d\theta}{d\nu}\right|\F_\tau\right]\right)\\
&\leq&\inf_{u\in
  \DD_\tau \cap L^2_a(\nu,H)}\EE_{\nu}\left[\left.f\circ
     W^u+\frac{1}{2}|u|_H^2\right|\F_\tau\right]\eeaa\nqed

\no Here is the main result.

\begin{theorem}
\label{3varrep}
Set $f:\W\rightarrow\RR$ measurable and $p,q>1$ such
that $p^{-1}+q^{-1}=1$ and $f\in L^p(\nu)$, $e^{-f}\in
L^q(\nu)$, then we have
$$-\log\EE_{\nu}\left[\left.e^{-f}\right|\F_\tau\right]=\inf_{u\in\D_\tau^i}\EE_{\nu}\left[\left.f\circ W^u+\frac{1}{2}|u|_H^2\right|\F_\tau\right]\;\;\nu-a.s.$$
\end{theorem}

\nproof The inequality
$$-\log\EE_{\nu}\left[\left.e^{-f}\right|\F_\tau\right]\leq\inf_{u\in
  \D_\tau^i}\EE_{\nu}\left[\left.f\circ
     W^u+\frac{1}{2}|u|_H^2\right|\F_\tau\right]$$
\no is an easy consequence of proposition \ref{3ineqvar}.
\no Let $\theta_0$ be the measure on W defined by
$$d\theta_0=\frac{e^{-f}}{\EE_{\nu}\left[\left.e^{-f}\right|\F_\tau\right]}d\nu$$
\no According to corollary \ref{3repr}, there exists $(u_n)\in
\left(\D_\tau^i\right)^\NN$ such that $\nu$-a.s.
\beaa \EE_{\nu}\left[\left.\frac{d W^{u_n}\nu}{d\nu}\log \frac{d W^{u_n}\nu}{d\nu}\right|\F_\tau\right] \rightarrow
 \EE_{\nu}\left[\left.\frac{d\theta_0}{d\nu}\log \frac{d\theta_0}{d\nu}\right|\F_\tau\right]\\
 \EE_{\nu}\left[\left.\frac{d W^{u_n}\nu}{d\nu}\log \frac{d\theta_0}{d\nu} \right|\F_\tau\right]\rightarrow
 \EE_{\nu}\left[\left.\frac{d\theta_0}{d\nu}\log \frac{d\theta_0}{d\nu}\right|\F_\tau\right]\eeaa
\no Denote $L_n=\frac{d W^{u_n}\nu}{d\nu}$, since $ W^{u_n}$
is $\nu$-a.s. invertible, we have
\beaa \EE_{\nu}\left[\left.f\circ W^{u_n}+\frac{1}{2}|u_n|_H^2\right|\F_\tau\right]=
\EE_{\nu}\left[L_n|\F_\tau\right]+\EE_{\nu}\left[L_n\log
 L_n|\F_\tau\right]\eeaa
\no When n goes to infinity, we have $\nu$-a.s.
$$\EE_{\nu}\left[L_n\log L_n|\F_\tau\right]\rightarrow
\EE_{\nu}\left[\left.\frac{d\theta_0}{d\nu}\log \frac{d\theta_0}{d\nu}\right|\F_\tau\right]$$
\no and since $f=-\log
\frac{d\theta_0}{d\nu}-\log\EE_{\nu}\left[\left.e^{-f}\right|\F_\tau\right]$, $\nu$-a.s.
\beaa\EE_{\nu}\left[fL_n|\F_\tau\right]&\rightarrow& \EE_{\nu}\left[\left.f\frac{d\theta_0}{d\nu}\right|\F_\tau\right]\eeaa
\no So finally, when n goes to infinity, $\nu$-a.s.
\beaa \EE_{\nu}\left[\left.f\circ
   W^{u_n}+\frac{1}{2}|u_n|_H^2\right|\F_\tau\right]&\rightarrow&\EE_{\theta_0}\left[f\right]+\EE_{\nu}\left[\left.\frac{d\theta_0}{d\nu}\log \frac{d\theta_0}{d\nu}\right|\F_\tau\right]
\\
&=&-\log\EE_{\nu}\left[\left.e^{-f}\right|\F_\tau\right]\eeaa
\no which conclude the proof.\nqed

\begin{theorem}
 Set $f:\W\rightarrow\RR$ a measurable function
verifying $\EE_{\nu}\left[|f|(1+e^{-f})\right]<\infty$, then if there
exists some $u\in \DD_\tau\cap L^2_a(\nu,H)$ such that $ W^u$ is
$\nu$-a.s. left-invertible and
$\frac{d W^u\nu}{d\nu}=\frac{e^{-f}}{\EE_{\nu}\left[\left.e^{-f}\right|\F_\tau\right]}$, then
we have
$$-\log\EE_{\nu}\left[\left.e^{-f}\right|\F_\tau\right]=\inf_{u\in \DD_\tau\cap L^2_a(\nu,H)}\EE_{\nu}\left[\left.f\circ W^u+\frac{1}{2}|u|_H^2\right|\F_\tau\right]\;\;\nu-a.s.$$
\end{theorem}

\nproof Since $ W^u$ is $\nu$-a.s. left invertible and that
$\frac{d W^u\nu}{d\nu}=\frac{e^{-f}}{\EE_{\nu}\left[\left.e^{-f}\right|\F_\tau\right]}$. We have
$$\frac{1}{2}\EE_{\nu}\left[\left.|u|_H^2\right|\F_\tau\right]=\EE_{\nu}\left[\left.\frac{e^{-f}}{\EE_{\nu}\left[\left.e^{-f}\right|\F_\tau\right]}\log
  \left(\frac{e^{-f}}{\EE_{\nu}\left[\left.e^{-f}\right|\F_\tau\right]}\right)\right|\F_\tau\right]$$
\no and
\beaa\EE_{\nu}\left[\left.f\circ
   W^u+\frac{1}{2}|u|_H^2\right|\F_\tau\right]&=&\EE_{\nu}\left[\left.\frac{e^{-f}}{\EE_{\nu}\left[\left.e^{-f}\right|\F_\tau\right]}f+\frac{e^{-f}}{\EE_{\nu}\left[\left.e^{-f}\right|\F_\tau\right]}\log\left(\frac{e^{-f}}{\EE_{\nu}\left[\left.e^{-f}\right|\F_\tau\right]}\right)\right|\F_\tau\right]\\
&=&-\log\EE_{\nu}\left[\left.e^{-f}\right|\F_\tau\right]\eeaa
\no and we conclude the proof with proposition \ref{3ineqvar}.\nqed

\begin{theorem}
Set $f:\W\rightarrow \RR$  a measurable function such that
$$-\log\EE_{\nu}\left[\left.e^{-f}\right|\F_\tau\right] = \inf_{u\in
  \DD_\tau\cap L^2_a(\nu,H)}\EE_{\nu}\left[\left.f\circ W^u+\frac{1}{2}|u|_H^2\right|\F_\tau\right]\;\;\nu-a.s.$$
\no Denote this infimum $J_*$. It is attained at $u\in \DD_\tau\cap L^2_a(\nu,H)$ if and only
if $ W^u$ is $\nu$-a.s. left-invertible and
$\frac{d W^u\nu}{d\nu}=\frac{e^{-f}}{\EE_{\nu}\left[\left.
      e^{-f}\right|\F_\tau\right]}$.
\end{theorem}

\nproof Denote $L=\frac{d W^u\nu}{d\nu}$. The direct implication is given by last theorem.
\no Conversely, if $ W^u$ is not $\nu$-a.s. left-invertible,
$\EE_{\nu}\left[L\log L|\F_\tau\right]<\frac{1}{2}\EE_{\nu}\left[\left.|u|_H^2\right|\F_\tau\right]$ and
\beaa -\log\EE_{\nu}\left[\left.e^{-f}\right|\F_\tau\right]&\leq& \inf_{\alpha\in
  \DD_\tau\cap L^2_a(\nu,H)}\EE_{\nu}\left[\left.f\circ W^\alpha+\frac{d W^\alpha\nu}{d\nu}\log\frac{d W^\alpha\nu}{d\nu}\right|\F_\tau\right]\\
&\leq&\EE_{\nu}\left[\left.f\circ W^u+L\log L\right|\F_\tau\right]\\
&<&\EE_{\nu}\left[\left.f\circ W^u+\frac{1}{2}|u|_H^2\right|\F_\tau\right]\eeaa
\no which is a contradiction.\newline
\no We get $L=\frac{e^{-f}}{\EE_{\nu}\left[\left.
      e^{-f}\right|\F_\tau\right]}$ by uniqueness of the minimizing
measure of $\inf_{\theta\in
  \mathcal{P}_\tau}\EE_\theta\left[\left.f+\log\frac{d\theta}{d\nu}\right|\F_\tau\right]$.\nqed

\section{\bf{Pr\'ekopa-Leindler theorem for conditional expectations}\label{pl3}}

\begin{definition}
We denote
$$H_b=\left\{h\in H, \dot{h}\;is\;dt-a.s.\;bounded\right\}$$
\end{definition}

\begin{remark}
Observe that $H_b\subset\D$ and that if $u\in\D$, $u(w)\in H_b$ $\nu$-a.s.
\end{remark}

\begin{theorem}
\label{3pl}
Assume that for any $u\in\D$,
$$W^{u}(w)=W^{u(w)}(w)\;\nu-a.s.$$
\no Set $t\in [0,1]$. Set $a,b,c:\W\rightarrow\RR$ positive and measurable such that for every $h,k\in H$ and $s\in [0,1]$ we have $\nu$-a.s.
$$a\circ  W^{ s h+(1- s)k}\exp\left(-\frac{1}{2}\left| s
  h+(1- s) k\right|_H^2\right)\geq\left(b\circ W^{h}\exp\left(-\frac{1}{2}\left|
  h\right|_H^2\right)\right)^{ s}\left(c\circ W^{k}\exp\left(-\frac{1}{2}\left|
  k\right|_H^2\right)\right)^{1- s}$$
\no then for any density d such that $h\in H_b \mapsto -\log d\circ  W^h$ is
$\nu$-a.s. concave and $\EE_{\nu}\left[d|\F_\tau\right]=1$, if $\theta$ denotes the measure on W given by
$\frac{dS}{d\nu}=d$, we have in $\bar{\RR}$:
$$\EE_\theta\left[a\left|\F_\tau\right.\right]\geq\left(\EE_\theta\left[b\left|\F_\tau\right.\right]\right)^ s\left(\EE_\theta\left[c\left|\F_\tau\right.\right]\right)^{1- s}$$
\end{theorem}

\nproof First observe that eventually replacing a,b,c with $da,db,dc$ and using Bayes formula
we only need to prove the case $d=1$ i.e. $\theta=\nu$\newline

\no With the convention $\log(\infty)=\infty$ and $\log(0)=-\infty$, we denote
$$\tilde{a}=-\log a,\tilde{b}=-\log b,\tilde{c}=-\log c$$
\no We begin with the case where there exists $m,M>0$ such that we
have $\nu$-a.s.
$$m\leq \tilde{a},\tilde{b},\tilde{c}\leq M$$
\no Set $ s\in [0,1]$ for $h,k\in H$, we have
\beaa &&a\circ  W^{ s h+(1- s)k}\exp\left(-\frac{1}{2}\left| s
  h+(1- s) k\right|_H^2\right)\\
  &&\geq\left(b\circ W^{h}\exp\left(-\frac{1}{2}\left|
  h\right|_H^2\right)\right)^{ s}\left(c\circ W^{k}\exp\left(-\frac{1}{2}\left|
  k\right|_H^2\right)\right)^{1- s}\eeaa
  \no So for $u_1,u_2\in\D_\tau^i$
  \beaa &&a\circ  W^{ s u_1+(1- s)u_2}\exp\left(-\frac{1}{2}\left| s
  u_1+(1- s) u_2\right|_H^2\right)\\
  &&\geq\left(b\circ W^{h}\exp\left(-\frac{1}{2}\left|
  h\right|_H^2\right)\right)^{ s}\left(c\circ W^{k}\exp\left(-\frac{1}{2}\left|
  k\right|_H^2\right)\right)^{1- s}\eeaa
\no hence applying the logarithm function, changing the sign and
taking the conditional expectation relative to $\F_\tau$ we obtain
\beaa &&\EE_{\nu}\left[\left.\tilde{a}\circ  W^{ s u_1 +(1- s)u_2}+\frac{1}{2}\left| s
  u_1+(1- s)u_2\right|^2_H\right|\F_\tau\right]\\
  &&\leq s\EE_{\nu}\left[\left.\tilde{b}\circ W^{u_1}+\frac{1}{2}\left|u_1\right|^2_H\right|\F_\tau\right]+(1- s)\EE_{\nu}\left[\left.\tilde{c}\circ W^{u_2}+\frac{1}{2}\left|u_2\right|^2_H\right|\F_\tau\right]\eeaa
\no So
\beaa \inf_{u\in \D_\tau^i}
&&\EE_{\nu}\left[\left.\tilde{a}\circ W^u+\frac{1}{2}|u|^2_H\right|\F_\tau\right]\\
&&\leq
   s\EE_{\nu}\left[\left.\tilde{b}\circ
     W^{u_1}+\frac{1}{2}\left|u_1\right|^2_H\right|\F_\tau\right]+(1- s)\EE_{\nu}\left[\left.\tilde{c}\circ
     W^{u_2}+\frac{1}{2}\left|u_2\right|^2_H\right|\F_\tau\right]\eeaa

\no According to theorem \ref{3varrep} we have
\beaa -\log\EE_{\nu}\left[\left.e^{-\tilde{a}}\right|\F_\tau\right]&\leq&
   s\EE_{\nu}\left[\left.\tilde{b}\circ
     W^{u_1}+\frac{1}{2}\left|u_1\right|^2_H\right|\F_\tau\right]+(1- s)\EE_{\nu}\left[\left.\tilde{c}\circ
     W^{u_2}+\frac{1}{2}\left|u_2\right|^2_H\right|\F_\tau\right]\eeaa
\no which implies
\beaa  -\log\EE_{\nu}\left[\left.e^{-\tilde{a}}\right|\F_\tau\right]&\leq&
   s\EE_{\nu}\left[\left.\tilde{b}\circ
     W^{u_1}+\frac{1}{2}\left|u_1\right|^2_H\right|\F_\tau\right]\\
    &&+(1- s)  \inf_{v\in \D_\tau^i} \EE_{\nu}\left[\left.\tilde{c}\circ
     W^{v}+\frac{1}{2}\left|v\right|^2_H\right|\F_\tau\right]\\
&=& s\EE_{\nu}\left[\left.\tilde{b}\circ
     W^{u_1}+\frac{1}{2}\left|u_1\right|^2_H\right|\F_\tau\right]-(1- s)\log\EE_{\nu}\left[\left.e^{-\tilde{c}}\right|\F_\tau\right]\eeaa
\no which implies once again
\beaa  -\log\EE_{\nu}\left[\left.e^{-\tilde{a}}\right|\F_\tau\right]&\leq& s\inf_{v\in \D_\tau^i}\EE_{\nu}\left[\left.\tilde{b}\circ
     W^{v}+\frac{1}{2}\left|v\right|^2_H\right|\F_\tau\right]\\
    &&-(1- s)\log\EE_{\nu}\left[\left.e^{-\tilde{c}}\right|\F_\tau\right]\\
&=&- s\log\EE_{\nu}\left[\left.e^{-\tilde{b}}\right|\F_\tau\right]
-(1- s)\log\EE_{\nu}\left[\left.e^{-\tilde{c}}\right|\F_\tau\right]\eeaa
\no taking the opposite and applying the exponential, we get
$$\EE_{\nu}\left[\left.e^{-\tilde{a}}\right|\F_\tau\right]\geq\left(\EE_{\nu}\left[\left.e^{-\tilde{b}}\right|\F_\tau\right]\right)^ s\left(\EE_{\nu}\left[\left.e^{-\tilde{c}}\right|\F_\tau\right]\right)^{1- s}$$

\no For the general case, denote for $n\in\NN$ and $m\in\NN^*$
\beaa &&\tilde{a}_n=\tilde{a}\wedge n,\tilde{b}_n=\tilde{b}\wedge n,\tilde{c}_n=\tilde{c}\wedge n\\
&&\tilde{a}_{nm}=\tilde{a}_n+\frac{1}{m},\tilde{b}_{nm}=\tilde{b}_n+\frac{1}{m},\tilde{c}_{nm}=\tilde{c}_n+\frac{1}{m}\eeaa
\no For every $h,k\in H$, we have $\nu$-a.s.:
$$\tilde{a}_{nm}\circ  W^{ s h +(1- s)k}+\frac{1}{2}\left| s
  h+(1- s)k\right|^2_H\leq  s\tilde{b}_{nm}\circ W^{h}+\frac{1}{2}\left|h\right|^2_H+(1- s)\tilde{c}_{nm}\circ W^{k}+\frac{1}{2}\left|k\right|^2_H$$
\no so the bounded case we treated above ensures that
$$\EE_{\nu}\left[\left.e^{-\tilde{a}_{nm}}\right|\F_\tau\right]\geq\left(\EE_{\nu}\left[\left.e^{-\tilde{b}_{nm}}\right|\F_\tau\right]\right)^ s\left(\EE_{\nu}\left[\left.e^{-\tilde{c}_{nm}}\right|\F_\tau\right]\right)^{1- s}$$

\no The monotone limit theorem enables us to take the limit with
relation to m and then to take it again with respect to n to get the
result.\nqed

\section{\bf{Examples}}

\no In this section we discuss several examples that fit into the framework we elaborated. Each time, we prove that the conditions of section \ref{fr3} and definition \ref{3DD} are satisfied, which ensure that every result from section \ref{fr3} to \ref{pl3} apply. We also discuss weather theorem \ref{3pl} applies or not.
\no See \cite{a1} for the omitted proofs concerning the diffusion, \cite{a2} for the omitted proofs concerning the other examples.

\subsection{\bf{Diffusion}}

\no Set $m\leq d\in \NN^*$ such that $m+d=n$, $c\in\RR$, $\sigma:\RR^m\rightarrow\M_{m,d}(\RR)$ and $b:\RR^m\rightarrow\RR^m$ bounded and lipschitz functions. $\sigma_i$ will denote the
i-th column of $\sigma$. Notice that every matrix
will be identified with its canonical linear operator. Set $(\Omega,\theta,(\mathcal{G}_t))$ a probability space, V a $\theta$-Brownian motion on
$\Omega$ with values in $\RR^d$. Set Y a $\RR^m$-valued strong solution of the
stochastic differential equation:
$$Y(t)=c+\int_0^t\sigma(Y(s))dV(s)+\int_0^tb(Y(s))ds$$
\no on $(\Omega,\theta,(\mathcal{G}_t),B)$. The hypothesis on $\sigma$ and b ensure the
existence and uniqueness of Y if we impose its paths to be continuous.\newline
\no We denote $\mu$
the Wiener measure on $C([0,1],\RR^d)$ and $\mu^X$ the measure on $C([0,1],\RR^m)$ the image measure of Y.
\newline\no We define the processes X and B on $\W$ by:
\beaa X(t)&:&(w,w')\in \W\mapsto w(t)\in\RR^m\\
B(t)&:&(w,w')\in \W\mapsto w'(t)\in\RR^d\eeaa

\begin{proposition}
Under $\mu^X\times\mu$, the law of X is $\mu^X$, B is a Brownian
motion and they are independent. There exists $\theta,\eta$ such that if we define $\beta_\grandX$ as
$$\beta_\grandX=\int_0^.\theta(X(s))dM(s)+\int_0^.\eta(X(s))dB(s)$$
\no $\beta_\grandX$ is a $\mu^X\times\mu$-Brownian motion and
$\mu^X\times\mu$-a.s.
$$X=c+\int_0^.\sigma(X(s))d\beta_\grandX(s)+\int_0^.b(X(s))ds$$
\end{proposition}

\no This construction of $\beta_\grandX$ is taken from \cite{ro}.

\begin{definition} We denote
$$\grandX=(X,\beta_\grandX)$$
\no and $\loigrandX$ its image measure.\newline
X is a $\loigrandX$ path-continuous strong solution of the
stochastic differential equation
\beaa X=&c+\int_0^.\sigma(X(s))d\beta_\grandX(s)+\int_0^.b(X(s))ds\eeaa
\no For $u\in G_0(\loigrandX,\beta_X)$, set $\beta_\grandX^u=\beta+u$ and $X^u$ the $\loigrandX$-a.s. path-continuous strong solution
of the stochastic differential equation
$$X^u=c+\int_0^.\sigma(X^u(s))d\beta_\grandX^u(s)+\int_0^.b(X^u(s))ds$$
\no Finally, we denote
\beaa\grandX^u&=&(X^u,\beta_\grandX+u)\eeaa
\end{definition}

\begin{theorem} $\left(\W,\loigrandX,\beta_\grandX,(\grandX^u)_{u\in\D}\right)$ verify the
    conditions of section \ref{fr3}. $\left(\W,\mu_a,\beta_\grandX,\left(\grandX^u\right)_{u\in G_0(\loigrandX,\beta_\grandX)}\right)$
    verify the conditions of definition \ref{3DD}.
\end{theorem}

\nproof (vii) of definition \ref{3DD} is clear, see \cite{a1} for the remainder of the proof.\nqed

\begin{corollary} It is clear that for every $u\in\D$, we clearly have $\loigrandX$-a.s.
$$\grandX^u(w)=\grandX^{u(w)}(w)$$
\no so theorem \ref{3pl} applies.
\end{corollary}

\subsection{\bf{Brownian bridge}}
\no We still denote $\mu$ the Wiener measure on $\W$. Set $a\in\R^n$, we denote $\mu_a$ the measure on $\W$ such that for any bounded measurable function f we have
$$\EE_{\mu_a}\left[f\right]=\EE_\mu\left[f|W_1=a\right]$$
\no $\mu_a$ can also be defined as follow : let $\mathcal{E}_a$ be the Dirac measure in a, $\mathcal{E}_a(W_1)$ is a positive Wiener distributions hence it defines a Radon measure $\nu_a$ on $\W$, then
$$\mu_a=\left(\frac{1}{2\pi}\right)^n\nu_a$$

\no We recall the definition of a Brownian bridge:
\begin{definition}
Set $(\Omega,\mathcal{G},Q)$ a probability space. An a-Brownian bridge X under a probability Q is a path-continuous Gaussian process such
that $\EE_Q\left[X(t)\right]=at$ and $cov(X(s),X(t))=\left((s\wedge
  t)-st\right)I_d$
\end{definition}

\begin{proposition}
W is an a-Brownian bridge under $\mu_a$, and the process $\beta_a$ defined as
$$\beta_a(t)=W(t)-at+\int_0^t\frac{W(s)-as}{1-s}ds$$
\no is a Brownian motion under $\mu_a$ and the filtrations of $\beta_a$ and
W completed with respect to $\mu_a$ are equal. Moreover, we have
$$W(t)=at+(1-t)\int_0^t\frac{d\beta_a(s)}{1-s}$$
\end{proposition}

\no The following remark will be useful in next section.\newline
\begin{remark}
For $a\in\RR^n$ and $t\in [0,1]$, we have $\mu_a$-a.s.
$$\beta_a(t)=W_t+\int_0^t\frac{W_s-a}{1-s}ds$$
\end{remark}

\begin{definition}
\no For $u\in G_0(\mu_a,\beta_a)$, we denote $\beta_a^u=\beta_a+u$.
\end{definition}

\begin{proposition}
\no Set $u\in G_0(\mu_a,\beta_a)$, then there exists a unique $\mu_a$-a.s. path
continuous process $W_a^u$ such that
$$W_a^u(t)=\beta_a^u(t)+at-\int_0^t\frac{W_a^u(s)}{1-s}ds$$
\no Furthermore, we have
\beaa W_a^u(t)&=&at+ (1-t)\int_0^t\frac{d\beta_a^u(s)}{1-s}\\
&=&W(t)+\int_0^t\left(\dot{u}(s)-\int_0^s\frac{\dot{u}(r)}{1-r}dr\right)ds\eeaa
\end{proposition}

\begin{theorem} $\left(\W,\mu_a,\beta_a,(W_a^u)_{u\in\D}\right)$ verify the
    conditions of section \ref{fr3}. $\left(\W,\mu_a,\beta_a,\left(W_a^u\right)_{u\in G_0(\mu_a,\beta_a)}\right)$
    verify the conditions of definition \ref{3DD}.
\end{theorem}
\nproof (vii) of definition \ref{3DD} is clear, see \cite{a2} for the remainder of the proof.\nqed

\begin{corollary} It is clear that for every $u\in\D$, we clearly have $\mu_a$-a.s.
$$W_a^u(w)=W_a^{u(w)}(w)$$
\no so theorem \ref{3pl} applies.
\end{corollary}

\subsection{\bf{Loop measure}}

\no We keep the notations of last section. Denote
$$S=\left\{a\in\RR^n,|a|=1\right\}$$
\no and set $\alpha:S \rightarrow \RR_+$ a locally lipschitz
function such that $\left\{x,\alpha(x)\neq 0\right\}$ is of strictly
positive measure for the Lebesgue measure on S and
$$\int_S\alpha(a)da=1$$
\no We define the measure $\nu_l$ as follow:
for any bounded measurable function f on W, we
set
$$\EE_{\nu_l}[f]=\int_S\alpha(a)\EE_{\mu_a}[f]da$$
\no For more on loop measures, see Fang's work in \cite{fan}.
\begin{definition}
We denote
\beaa h_a:(t,x)&\in&
[0,1)\times\RR^n\mapsto\left(\frac{1}{\pi(1-t)}\right)^{\frac{n}{2}}\exp\left(\frac{-\left|x-a\right|^2}{2(1-t)}\right)\\
h:(t,x)&\in& [0,1)\times\RR^n\mapsto\int_S \alpha(a)h_a(t,x)da\eeaa
\end{definition}

\begin{proposition}
\label{3dl}
Set $a\in\RR^n$ and $t\in[0,1)$, then
$$\left.\frac{d\mu_a}{d\mu}\right|_{\F^W_t}=h_a(t,W_t)$$
\end{proposition}

\begin{proposition}
Set $t\in [0,1)$, we have
$$\left.\frac{d\nu}{d\mu}\right|_{\F_t^W}=h(t,W_t)$$
\end{proposition}

\begin{proposition}
Define
$$\beta_{\lo}(t)= W(t)-\int_0^t\frac{h'(s,W(s))}{h(s,W(s))}ds$$
\no where h' designates the partial derivative of h with respect to
x.\newline
\no Then $\beta_{\lo}$ is a $\nu_l$ Brownian motion
and the filtrations of W and $\beta_{\lo}$ completed with respect to
$\nu_l$ are equal.
\end{proposition}

\begin{definition}
\no For $u\in G_0(\nu_l,\beta_{lo})$, we denote $\beta_{\lo}^u=\beta_{\lo}+u$.
\end{definition}

\begin{proposition}
\no Set $u\in G_0(\nu_l,\beta_{lo})$, then there exists a unique $\nu_l$-a.s. path
continuous process $W_{\lo}^u$ such that
$$W_{\lo}^u(t)=\beta_{\lo}^u(t)+\int_0^t\frac{h'(s,W_{lo}^u(s))}{h(s,W_{lo}^u(s))}ds$$
\end{proposition}

\begin{theorem} $\left(\W,\nu_l,\beta_{\lo},(W_{\lo}^u)_{u\in\D}\right)$ verify the
    conditions of section \ref{fr3}. $\left(\W,\nu_l,\beta_{lo},\left(W_{\lo}^u\right)_{u\in G_0(\nu_l,\beta_{lo})}\right)$
    verify the conditions of definition \ref{3DD}.
\end{theorem}
\nproof (vii) of definition \ref{3DD} is clear, see \cite{a2} for the remainder of the proof.\nqed

\begin{corollary} It is clear that for every $u\in\D$, we clearly have $\nu_l$-a.s.
$$W_{lo}^u(w)=W_{lo}^{u(w)}(w)$$
\no so theorem \ref{3pl} applies.
\end{corollary}

\subsection{\bf{Diffusing particles without collision}}

\no Set $\sigma,b,\delta,\gamma\in\RR$ such that
$$\sigma^2\leq 2\gamma$$
\no The proof of the following theorem can be found in \cite{shi} or \cite{cl}.
\begin{theorem}
\label{3par}
Set $(\Omega,\theta,(\mathcal{G}_t))$ a filtered probability space, $(z_1(0),...,z_n(0))\in\RR^n$ and $B=(B_1,...,B_n)$ a $\RR^n$-valued $\theta$-Brownian motion. We consider the following stochastic differential system:
\beaa Z_1(t)&=&z_1(0)+\sigma B_1(t)+b\int_0^t Z_1(s)ds+ct+\gamma\sum_{j\in\{1,...,n\}\backslash\{1\}}\int_0^t\frac{ds}{Z_1(s)-Z_j(s)}\\
\vdots&&\\
Z_n(t)&=&z_n(0)+\sigma B_n(t)+b \int_0^tZ_n(s)ds+ct+\gamma\sum_{j\in\{1,...,n\}\backslash\{n\}}\int_0^t\frac{ds}{Z_n(s)-Z_j(s)}\eeaa
\no under the condition that $\theta$-a.s. for every $t\in [0,\infty)$
$$Z_1(t)\leq...\leq Z_n(1)$$
\no This system admits a unique strong solution on $(\Omega,\theta,(\mathcal{G_t}),B)$ and the first collision time is $\theta$-a.s. equal to $\infty$.
\end{theorem}

\no Consider $(\Omega,\theta,(\mathcal{G}_t))$ a filtered probability space, $(z_1(0),...,z_n(0))\in\RR^n$ and $B=(B_1,...,B_n)$ a $\RR^n$-valued $\theta$-Brownian motion, and Z the strong solution of the stochastic differential system of theorem \ref{3par}. Denote $\nu_{pa}=Z$ the image measure of Z. For $1\leq i\leq n$, denote $W_1,...,W_n$ the coordinates of W and define

$$M_i(t)=W_i(t)-z_i(0)-b\int_0^tW_i(s)ds-ct-\gamma\sum_{j\in\{1,...,n\}\backslash\{i\}}\int_0^t\frac{ds}{W_i(s)-W_j(s)}$$
\no and
$$M=(M_1,...,M_n)$$
\no M is a local martingale and
$$\langle M_i,M_j\rangle(t)=\sigma^2t$$
\no Define
$$\beta_{pa}=\frac{1}{\sigma}M$$
\no Levy theorem clearly ensures that $\beta$ is a $\nu_{pa}$-Brownian motion and we clearly have for every $1\leq i\leq n$,
$$W_i(t)=z_i(0)+\sigma\beta_{pa,i}(t)+b\int_0^tW_i(s)ds+ct+\gamma\sum_{j\in\{1,...,n\}\backslash\{i\}}\int_0^t\frac{ds}{W_i(s)-W_j(s)}$$

\no For $u\in G_0(\nu_{pa},\beta_{pa})$ denote
$$\beta_{pa}^u=\beta_{pa}+u$$
\no and $\nu_{pa}^u$ the probability measure given by
$$\frac{d\nu_{pa}^u}{d\nu_{pa}}=\rho(-\delta_{\beta_pa}u)$$
\no According to Girsanov theorem, $\beta_{pa}+u$ is a Brownian motion under $\nu_{pa}^u$, so according to theorem \ref{3par}, there exists a unique $\nu_{pa}^u$-a.s. continuous process $W_{pa}^u=(W_{pa,i}^u,...,W_{pa,n}^u)$ such that $\nu_{pa}^u$-a.s. for every $1\leq i\leq n$
$$W_{pa,i}^u(t)=z_i(0)+\sigma\beta_{pa,i}^u(t)+b\int_0^tW_{pa,i}^u(s)ds+ct+\gamma\sum_{j\in\{1,...,n\}\backslash\{i\}}\int_0^t\frac{ds}{W_{pa,i}^u(s)-W_{pa,j}^u(s)}$$
\no and $\nu_{pa}^u$-a.s. for every $t\in[0,1]$
$$W_{pa,1}^u(t)\leq...\leq W_{pa,n}^u(t)$$
\no Since $\nu_{pa}^u\sim\nu_{pa}$, $W^u$ is $\nu_{pa}$-a.s. continuous and $\nu_{pa}$-a.s. for every $1\leq i\leq n$
$$W_{pa,i}^u(t)=z_i(0)+\sigma\beta_{pa,i}^u(t)+b\int_0^tW_{pa,i}^u(s)ds+ct+\gamma\sum_{j\in\{1,...,n\}\backslash\{i\}}\int_0^t\frac{ds}{W_{pa,i}^u(s)-W_{pa,j}^u(s)}$$
\no and $\nu_{pa}$-a.s. for every $t\in[0,1]$
$$W_{pa,1}^u(t)\leq...\leq W_{pa,n}^u(t)$$

\begin{theorem} $\left(\W,\nu_{pa},\beta_{\pa},(W_{\pa}^u)_{u\in\D}\right)$ verify the
    conditions of section \ref{fr3}.\newline \no  $\left(\W,\nu_{pa},\beta_{pa},\left(W_{\pa}^u\right)_{u\in G_0(\nu_{pa},\beta_{pa})}\right)$
    verify the conditions of definition \ref{3DD}.
\end{theorem}

\nproof (vii) of definition \ref{3DD} is clear, see \cite{a2} for the remainder of the proof.\nqed

\begin{corollary} It is clear that for every $u\in\D$, we clearly have $\nu_{pa}$-a.s.
$$W_{pa}^u(w)=W_{pa}^{u(w)}(w)$$
\no so theorem \ref{3pl} applies.
\end{corollary}

\vspace{2cm}
\footnotesize{
\noindent
K\'evin HARTMANN, Institut Telecom, Telecom ParisTech, LTCI CNRS D\'ept. Infres, \\
23 avenue d'Italie, 75013, Paris, France\\
kevin.hartmann@polytechnique.org}

\end{document}